\documentclass[11pt]{amsart}
\usepackage{fullpage,amsfonts,amsmath,amscd,amssymb}
\theoremstyle{plain}

\newtheorem*{thm}{Theorem}

\theoremstyle{remark}

\newtheorem*{rems}{Remarks}

\DeclareMathOperator{\red}{red}
\newcommand{\ep}{\epsilon}
\newcommand{\C}{\mathbb{C}}

\newcommand{\EO}{\mathcal{O}_{\epsilon}[G]}   

\title{$\EO$ is a free module over $\mathcal{O}[G]$}
\author{Kenneth A. Brown}
\address{Department of Mathematics, University of Glasgow, 
Glasgow G12 8QW}
\email{kab@maths.gla.ac.uk}

\author{Iain Gordon}
\address{Department of Mathematics, University of Glasgow, 
Glasgow G12 8QW}
\email{igordon@Mathematik.Uni-Bielefeld.DE}

\author{J. T. Stafford} 
\address{Department of Mathematics, University of Michigan,
Ann Arbor, MI 48109-1109, USA.}
\email{jts@math.lsa.umich.edu}

\thanks{All three authors were visiting and partially supported by MSRI while
this research was conducted and extend their thanks to that institution.
The research for this paper was undertaken 
 while the second author was 
 supported by TMR grant ERB FMRX-CT97-0100 at the 
 University of Bielefeld and the third author supported 
 in part by the NSF}

\begin{document}
\begin{abstract}
We show that the quantised function algebra $\EO$ of a simply 
connected semisimple algebraic group $G$ at a root of unity is a free 
module over the subring isomorphic to $\mathcal{O}[G]$.
\end{abstract}
\maketitle 

Let $G$ be a simply-connected semisimple
algebraic group over $\C$. Let $\ell > 1$ be an odd integer, prime
to $3$ if $G$ has a component of type $G_2$, and let $\ep\in \C$
be a primitive $\ell^{\text{th}}$ root of unity. The quantised
function algebra of $G$ at $\ep$, denoted $\EO$, is a noetherian
$\mathbb{C}$-algebra containing the ring of regular functions of
$G$, denoted $\mathcal{O}[G]$, in its centre, \cite{declyu1}.
Since a false proof of the following theorem, and a proof of the
special case $G = SL(2)$, have both recently appeared in the
literature (see the remarks below for details), it seems
worthwhile to record the full result.
\begin{thm}
As a module over $\mathcal{O}[G]$, the algebra 
$\EO$ is free of rank $\ell^{\dim G}$.
\end{thm}
\begin{proof}
Thanks to \cite[Theorem~7.2]{declyu1} $\EO$ is
 a projective $\mathcal{O}[G]$-module of rank 
 $\ell^{\dim G}$. By a result of Marlin, \cite[Corollaire~3]{marlin}, 
 the Grothendieck group of projective modules over
  $\mathcal{O}[G]$ is trivial, in other words
\[
K_0(\mathcal{O}[G]) \cong \mathbb{Z}.
\]
In particular, if $P$ is a finitely generated projective
$\mathcal{O}[G]$-module then $P$ is stably free.  Hence if the
rank of $P$ is greater than the Krull dimension of
$\mathcal{O}[G]$, then $P$ is necessarily free, 
\cite[Corollary~IV.3.5]{Bass}. Since $\ell
>1$ we have $\text{Kdim}\mathcal{O}[G] = \dim G < \ell^{dim G} =
\text{rank}\EO$, so the theorem follows.
\end{proof}

It is incorrectly stated in \cite[Lemma~8]{rum} that this result follows from
\cite[Theorem~2.2]{sch} in the more general setting of a Hopf algebra, $U$, finitely
generated over a central sub-Hopf algebra, $O$. However, there exist numerous
examples in the literature of Hopf algebras that are not free over central
sub-Hopf algebras. The ones closest in spirit to the present work occur in
\cite{wat}, where the author shows that, when $n$ is even, $U=\mathcal{O}[SL_n(\mathbb C)]$ is not
free over the subring  $O=\mathcal{O}[PSL_n]$.  For example, consider the case $n=2$. Then, $O$ is the fixed ring
$U^A$, where $A=\mathbb Z/2\mathbb Z=\langle \sigma \rangle$ acts on the
generators $x_{ij}$ of $U$ by $\sigma(x_{ij})=-x_{ij}$. This is even a
Hopf-Galois extension with central invariants for the Hopf algebra $H=\mathbb
CA$, in the notation of \cite[\S1.1]{rum}. This requires that $U$ is an
$H$-comodule algebra (use the map  $\rho : U\to U\otimes H$ defined by $x_{ij}
\mapsto x_{ij}\otimes \sigma$) such that $U^H = \{x\in U : \rho(x)=x\otimes 1\}
= O$,  that the natural map $U\otimes_OU \to U\otimes H$ given by  $x\otimes
y\mapsto (x\otimes 1)\rho(y)$ is bijective, together with certain naturality 
conditions. These are easy to check in this case.
\begin{rems}
1. It is clear the proof of the theorem generalises a little. Assume $k$ is an algebraically closed field and $U$ is a noetherian prime Hopf $k$-algebra, finitely generated as a module over the central sub-Hopf algebra $O$. Then $U$ is a projective $O$-module, \cite[Theorem~1.7]{kretak}. If $K_0(O) \cong \mathbb{Z}$ then freeness follows as above when $\dim O$ is less than the rank of $U$ over $O$.
\\
2. One way to check that $K_0(O) \cong \mathbb{Z}$ is as follows. The algebra $O$ is the ring of regular functions of an irreducible affine algebraic group, say $G$,  \cite[Section~I.3]{hoc}. Let $R_u(G)$ be the unipotent radical of $G$ and set $G_{\red} = G/R_u(G)$, by definition a reductive group. Thanks to \cite[Proposition~4.1]{quil}, the projection $G\longrightarrow G_{\red}$ induces an isomorphism in $K$-theory, $K_0(G) \cong K_0(G_{\red})$. Now, if the commutator subgroup of $G_{\red}$, a semisimple algebraic group, is simply-connected we have an isomorphism $K_0(G_{\red}) \cong \mathbb{Z}$, \cite[Corollary~1.7 and Corollary~4.7]{mer}. 
\\
3. For other situations where 
 a Hopf algebra is free over particular subalgebras, see for example, 
 \cite{Ra}.
\\ 4.
\cite{iti} proves the theorem for the case $G=SL(2)$; in this case
an explicit free basis is provided.
\\ 5. We do not know whether $\EO$ is a cleft extension of $\mathcal{O}[G]$. It would be interesting to find general conditions implying this.
\\ 6. The theorem appears in \cite{brogor2} too, where it is used in studying the representation theory of quantised function algebras at roots of unity.
\end{rems}


\begin{thebibliography}{10}

\bibitem{brogor2}
K.A. Brown and I.~Gordon,
\newblock The ramification of centres: quantised function algebras at roots of
  unity,
\newblock math.RT/9912042.

\bibitem{Bass}
H. Bass,
\newblock {\em Algebraic K-theory},
\newblock Benjamin, New York, 1968.

\bibitem{iti}
L. Dabrowski, C. Reina and A. Zampa,
\newblock $A[Sl_q(2)]$ at roots of unity is a free module over $A[Sl(2)]$,
\newblock  math.QA/0004092.

\bibitem{declyu1}
C.~De~Concini and V.~Lyubashenko,
\newblock Quantum function algebras at roots of 1,
\newblock {\em Adv. Math.}, {\bf 108} (1994), 205--262.

\bibitem{hoc}
G.P. Hochschild, {\em Basic theory of algebraic groups and Lie algebras}, { Graduate Texts in Mathematics} 75, Springer-Verlag, Berlin, 1981.

\bibitem{kretak}
H.F. Kreimer and M. Takeuchi, Hopf algebras and Galois extensions of an algebra, {\em Indiana U. Math. J.}, {\bf 30} (1981), 675--692.

\bibitem{marlin}
R. Marlin, Anneaux de Grothendieck des vari\'{e}t\'{e}s de drapeaux, {\em Bull. Soc. math. France}, {\bf 104} (1976), 337-348.

\bibitem{McC-Rob}
J.C. McConnell and J.C. Robson,
\newblock {\em Noncommutative Noetherian Rings},
\newblock J. Wiley and Sons, 1988.

\bibitem{mer}
A.S. Merkur'ev,
\newblock Comparison of the equivariant and the ordinary $K$-theory of algebraic varieties, {\em St. Petersburg Math. J.} {\bf 9} (1998), 815--850.

\bibitem{quil}
D. Quillen,
\newblock Higher algebraic $K$-theory. I, {\em Algebraic $K$-theory $I$: Higher $K$-theories}, Springer Lecture Notes in Mathematics 341, (1973), 85--147. 

\bibitem{Ra}
D.~E.~Radford, Freeness (projectivity) criteria for Hopf algebras 
over Hopf subalgebras, {\em J. Pure Appl. Algebra}, {\bf 11}
(1977/78),  15--28.

\bibitem{rum}
D.~Rumynin,
\newblock Hopf-{G}alois extensions with central invariants and their geometric
  properties,
\newblock {\em Algebra and Rep. Theory}, {\bf 1} (1998), 353--381.

\bibitem{sch}
H-J. Schneider,
\newblock Normal basis and transitivity of crossed products
for Hopf algebras,
\newblock {\em J. Algebra}, {\bf 152} (1992), 289--312.

\bibitem{wat}
W.~C.~Waterhouse, The module structure of certain Hopf algebra extensions,
{\em Comm. Algebra}, {\bf 10} (1982), 115-120.
\end{thebibliography}
\end{document}